\documentclass[preprint,12pt]{elsarticle}

\usepackage{algpseudocode}
\usepackage{algorithm}

\usepackage{graphicx}

\usepackage{amssymb}

\usepackage[cmex10]{amsmath}

\journal{Nuclear Instruments and Methods B}

\begin{document}

\begin{frontmatter}



\title{PyHST2: an hybrid distributed code for
 high speed tomographic reconstruction with iterative reconstruction 
and a priori knowledge capabilities}
\author[1]{Alessandro Mirone\corref{cor1}}
\author[2]{Emmanuelle Gouillart}
\author[1,3]{Emmanuel Brun}
\author[1]{Paul Tafforeau}
\author[1]{Jerome Kieffer}
\cortext[cor1]{mirone@esrf.fr}

\address[1]{ESRF 6 Rue Jules Horowitz 38000 Grenoble, France}
\address[2]{Joint Unit CNRS/Saint-Gobain Surface du Verre et Interfaces 39 quai Lucien Lefranc 93303 Aubervilliers, France}
\address[3]{Department of Physics, Ludwig Maximilians University, Am Coulombwall 185748 Garching, Germany}


\begin{abstract}
We present the PyHST2 code which is in service at ESRF for phase-contrast and absorption tomography. This code has been engineered to sustain the high data flow typical of the $3^{rd}$ generation synchrotron facilities (10 terabytes per experiment) by adopting a distributed and pipelined architecture. The code implements, beside a default filtered backprojection reconstruction, iterative reconstruction techniques with {\em a-priori} knowledge. These latter are used to improve the reconstruction quality or in order to reduce the required data volume and reach a given quality goal. The implemented   {\em a-priori} knowledge techniques are based on  the total variation penalisation and
a new recently found convex functional which is based on overlapping patches. 
 We give details of the different methods and their implementations while the code is distributed under free license. 
 We provide methods for estimating, in the absence of  ground-truth data, the optimal parameters values for {\em a-priori} techniques.
\end{abstract}


\end{frontmatter}


\section{introduction}
\label{introduction}
Progresses in both digital detector technology and in the X-ray beam lead to increasing data rates produced by modern synchrotron beamlines. 
In the case of X-ray micro-tomography,  data sets of some tens of terabytes are routinely acquired for an experiment of few days~\cite{Brunet2004277}. 
In order to analyze the produced data, special solutions must be adopted to fit the limits of available input-output bandwidth and computing power. On another hand, users of X-ray tomography aim to push even further the frontiers of their studies towards new domains which require finer time resolution~\cite{Salvo2003}, better signal to noise ratio, and less radiation damage as in the case of medical tomography~\cite{Zhao2012}. 
While waiting for further progresses in detection technology, these requirements must be satisfied as much as possible with available tools.


We present  the PyHST2 code which solves the High Performance Computing (HPC) issues with a distributed hybrid architecture and applies  {\em a-priori} knowledge
techniques to obtain high quality reconstructions even with reduced amount of data and in the presence of noise.
The details of the software, and  its HPC implementation, are exposed in section \ref{software}.


In section \ref{apriori} we present  two {\em a-priori} knowledge techniques that are available in PyHST2. These two techniques are the total variation regularization, on one side, which is well suited to piece-wise constant samples and a new overlapping-patches technique~\cite{overlap}, on the other side, that is suited to samples for which a sparse representation can be learned by the dictionary-learning technique.
Some applications of the implemented  state-of-the-art techniques are illustrated in section \ref{application} . 
 



\section{Software Description}
\label{software}
\subsection{Experimental scope}

The PyHST2 code has been conceived for synchrotron beamline setups and relies on the  assumptions that the beam is monochromatic, has limited divergence and that we can consider kinematic propagation through the sample.

The reconstruction algorithm needs a set of projection spanning a sample-rotation angle of at least 180 degrees.  When the sample rotation spans 360 degrees  the effective diameter of the detector can be  virtually doubled, as an option, by the reconstruction algorithm if, in the experimental set-up,  the rotation axis projection is close to the detector border.       
 The PyHST code performs reconstruction from raw data. These data comprise, beside the radiographies, the beam profile that is obtained by taking  radiographies in the absence of the sample and  the detector noise which is function of the detector pixel. In case of mechanical defects, the rotation axis displacement and the angular position can be given in input, as a function of  the projection, to correct the reconstruction. 
  In the case of propagation enhanced phase contrast tomography we consider that the beam is spatially coherent. We consider also, for this case, that the effect of the beam non-uniformity is small enough so that it can be factored out from the radiography  by a simple division. This means neglecting the sample-detector propagation effects that a strong non-uniformity could have on the signal from the sample.

Concerning the beam divergence we are currently adapting the code to nano-tomography experiments“\cite{nano} where a finite, thought small,  divergence is assumed. This application will be available in
future versions while the currently distributed version of PyHST still implements parallel geometry.

The kinematical approach hypothesis is valid when the resolution  (the linear size $d$ of a volume voxels) is far from the diffraction limit : $ d^2 \gg \lambda L $  given by the wavelength $\lambda$ and the sample diameter $L$ .

The sample-to-detector distance can instead be large enough to enhance the phase-contrast. We apply in this case the Paganin filter~\cite{paganin}.
The experimental situation where a limited amount of data is available, or where noise must be reduced, can be treated with the {\em a-priori} knowledge technique exposed in section \ref{apriori}.

\subsection{Implementation Details}
\label{Implementation Details}

The fundamental hardware unit that we consider for the deployment of the calculation is a multi-core processor (CPU) with its optional 
graphic card accelerator (GPU).
  The PyHST2 code spawns one process per  CPU. Each process takes full control of the CPU cores and of its associated GPU if present. Every process is bound to its associated CPU manipulating its CPU affinity by means of the Linux  {\it taskset} command.  
The controlling process is a script written in Python while the controlled threads are written in C language, for the CPUs, and in CUDA for the GPUs.
The controlling script takes care of optimally splitting the calculation in order to use the maximum amount of memory available avoiding swapping on virtual memory, and of always keeping busy the processing resources.
 The total data volume is often larger than the available memory. Therefore the calculation is done step by step, reading at each step only a part of 
the data, and reconstructing the sample sub volume by sub volume. The convolution kernel used for phase contrast treatment has often a large radius. For this reason the read data volume, for a single step, is larger than the reconstructed  sub volume.
 The communication between threads is direct, from thread to thread, without passing 
by the Python level. This is done either via the shared memory if the threads belong to the same process or via MPI communications if they are on different processes or  hosts. The communications between units is used to transpose the treated data into sinograms.
 It will also  be used in planned future developments for 3D regularisations.

 No number crunching or 
data transfer is done by the  Python routines.  These routines  control the sequence of data packets using  an attributed numerical identifier and using pointers to the available  memory regions.
   Each controlling process accounts also for  all the computing resources (cores and GPUs) associated to the CPU the process is running on.   All the resources are organized at 
the {\bf C} level into a global hierarchical {\bf C } structure which is initialized from Python. When the resources necessary for a given  pending computing task become available, the Python controlling process calls the corresponding  processing object. This consists in calling a  {\bf C} routine  which takes  as arguments a pointer to the global hierarchical  {\bf C} structure and the identifiers of the resources to be used. The {\bf C} routines are wrapped using the Python {\bf C} Application Program Interface  (C-API). The Python Global Interpreter Lock (GIL) is released by the called routines, thus realizing true multi-threading.

\subsection{Processing Pipe}

The data undergo several preprocessing steps before application of the reconstruction algorithm.
The first treatments consist in the removal of spurious signal coming from the detector read-out dark-current background-noise and in the correction for the beam spatial non-homogeneity, called flat-field correction. The dark-current is an image giving the averaged  signal measured by the detector when no beam is present. 
 The flat-field images are recorded in the absence of the sample: they are  images of the beam. During the tomographic scan several flat-fields can be acquired at different times to track the beam shape drifts.
The flat field corrections are applied dividing every radiography by the flat-field obtained by linear interpolation between the different flat-fields  acquisition times.   
 After correction for the flat field the signal can still contain spurious contributions. One of these is due to hot-spots in the detector. In the case of CCD detectors certain pixels may produce a constant noise.  The contribution
of these pixels is regularized by a median filter:  if the median of the pixel neighbors is smaller than the pixel value by more than a given threshold, then the pixel value is substituted with the neighbors median.

   When the sample-detector distance is long enough, the phase contrast becomes visible and can improve the signal/noise ratio by several order of magnitudes in the cases where the absorption contrast tends to zero. In the case of samples having an homogenous $\delta/\beta$ ratio (where the optical index is $n=1-\delta+i \beta$) both the amplitude and phase transmitted through the sample can be retrieved from the radiography. This is realized by implementation of the Paganin method~\cite{paganin} and  consists in the application of a two-dimensional convolution. 

After preprocessing of the radiographies, the resulting data volume, which is  initially a sequence of radiographies,  is rearranged  as a sequence of sinograms :  one sinogram per  detector line perpendicular to the rotation axis. This rearrangement is particularly convenient in the parallel geometry case,  where a given reconstructed line depends on one and only one sinogram.

\begin{figure}
\caption{\label{rings} Rings correction on a mouse leg tomography slice : a) with no corrections residual pixel-based errors give concentrical ring artefacts.  b) the application of the rings correction filter removes the rings but new 
introduces artefacts which are due to the bones. c) a pro-filter thresholding treatement of the bone suppresses the new artefacts. }

\begin{minipage}[t]{0.33\textwidth}%
\begin{center}
\includegraphics[width=0.95\columnwidth,height=0.95\columnwidth]{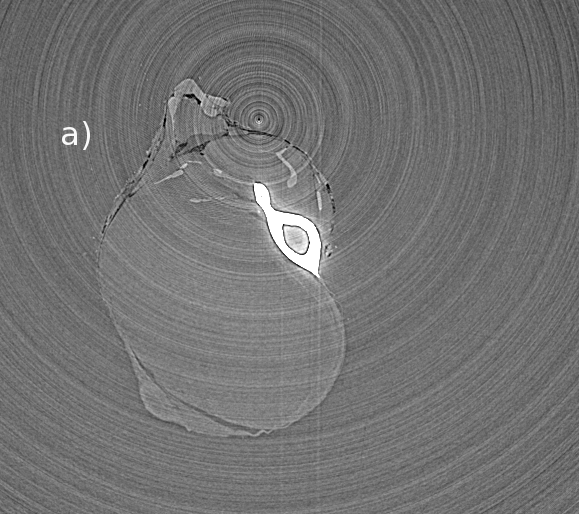}
\par\end{center}%
\end{minipage}
\hfill{}%
\begin{minipage}[t]{0.33\textwidth}%
\begin{center}
\includegraphics[width=0.95\columnwidth,height=0.95\columnwidth]{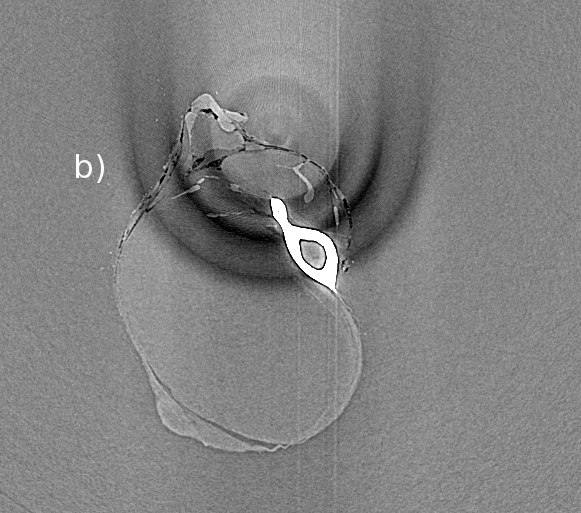}
\par\end{center}%
\end{minipage}\hfill{}%
\begin{minipage}[t]{0.33\textwidth}%
\begin{center}
\includegraphics[width=0.95\columnwidth,height=0.95\columnwidth]{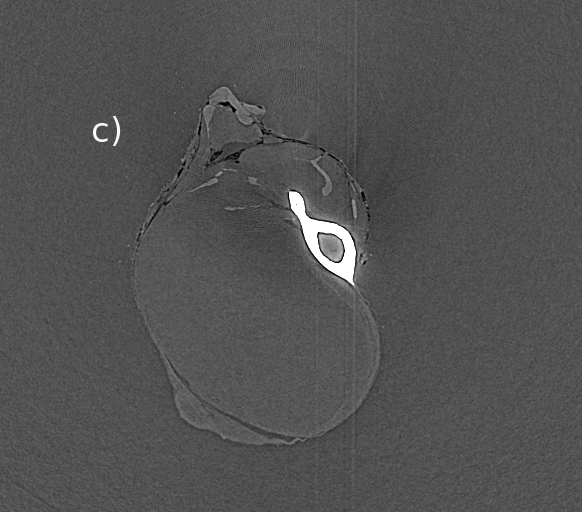}
\par\end{center}%
\end{minipage}
\end{figure}

The last preprocessing step is performed on the sinograms to reduce the so-called ring artifacts.
These artifacts are due to residual errors which are linked to a given pixel that still remain in the treated data even after applications of flat-field corrections and hot-spots removal filter. 
These errors are visible in the sinogram as features which are  parallel to the projection-angle axis and give raise, in the reconstructed slice, to ring-shaped artefacts.  
An optional filter can be used for these features. The filter reduces first the sinogram to a 1D signal by summing over the projection-angles. An high-frequencies filter is then applied in an attempt to extract these features and the result is subtracted from the sinogram. This approach can  give satisfactory results but, in some cases,
new artefacts can appear if the filtered signal contains a part coming from the sample.
We show in figure~\ref{rings}a(left) a tomographic slice of a mouse leg with no corrections for the ring artefacts.
We have applied the ring corrections for figure~\ref{rings}b(center). This corrects well the ring artefacts
but  new artefacts appear which are linked to the sample, not to the detector. In particular to the leg bone which has a strong absorption
and sharp borders.  PyHST2 can treat this problem by thresholding the bone in a first reconstructed slice. The thresholded signal is then subtracted to the data before application of the high-pass filter. The 
result is shown in figure~\ref{rings}c(right).

The reconstruction is done either by filtered back projection or by the advanced {\em a-priori} knowledge technique described in the next section.
The back projection has an efficient GPU implementation storing the sinogram into a GPU texture
while the  reconstructed volume is decomposed into tiles. Each  tile is reconstructed by a block
of threads. The values to be back-projected  are read from the texture using GPU hardware
interpolation. The tiles decomposition ensures  an high cache hit ratio thanks to the texture cache 2D locality.
The forward projection, used in the iterations for the {\em a-priori} algorithms, is done instead
keeping the volume on a texture.  


\section{A-priori Knowledge Techniques}
\label{apriori}

\subsection{Total Variation Penalization}
\label{aprioritv}

Classical filtered back-projection does not integrate any {\em a-priori}
information on the absorption image $\mathbf{x}$ during the
reconstruction. For noisy or under-sampled data, this has the consequence
that the projection  of the image $\mathbf{x}$ on eigenvectors of
$\mathbf{P}^\mathsf{T}\mathbf{P}$
of null or small eigenvalue (where $\mathbf{P}^\mathsf{T}$ is the back-projection operator)  is not
correctly reconstructed. Such indeterminacy arises for example during
in-situ imaging of evolving systems~\cite{Salvo2003}, where one strives
to reduce the acquisition time at the expense of either the signal to
noise ratio of the radiographies, or of the number of radiographies. This
difficulty can be overcome by using information on the image, such as the
sparsity of the image (e.g. for very porous materials) or its spatial
regularity. In a general setting, the reconstruction problem then amounts
to an optimization problem:
\begin{equation}
    \mathbf{x} = \mathrm{argmin}_{\mathbf{x}}\, \frac{1}{2}||\mathbf{y} -
    \mathbf{P}\mathbf{x}||^2 + f(\mathbf{x}),
\end{equation}
where $f$ incorporates the information on the image and $\mathbf{y}$ are the data. In the 2000's,
supplementing missing measures by
incorporating prior information has won its mathematical spurs within the
field of \emph{compressed sensing}~\cite{Candes2006,Candes2006b}. For
materials with a few phases of constant absorption, a classical
penalization function $f$ in tomographic reconstruction is the total
variation of the image~\cite{Song2007, Sidky2008, Herman2008,
Tang2009, Sidky2010, Jia2010}, that is a
$\ell_1$ norm of the image gradient. The use of the $\ell_1$ norm
promotes gradient sparsity, hence piecewise-constant images. To avoid
staircasing effects, we use here the isotropic total
variation~\cite{Chambolle2004} that promotes joint sparsity of all components of the gradient
together:
\begin{equation}
    TV(\mathbf{x}) = \sum_{\mathrm{pixels}}
    \sqrt{(\partial_1\mathbf{x})^2 + (\partial_2\mathbf{x})^2}.
\end{equation}

The optimization problem solved is then
\begin{equation}
        \mathbf{x} = \mathrm{argmin}_{\mathbf{x}}\,
	\frac{1}{2}||\mathbf{y} -
	\mathbf{P}\mathbf{x}||^2 + \beta TV(\mathbf{x}).
\label{eq:opt_tv}
\end{equation}
The parameter $\beta$ controls the relative importance of the
data-fidelity term and the spatial regularization: the higher $\beta$,
the smoother the reconstructed image. The optimization problem
(\ref{eq:opt_tv}) is convex, but non-smooth because of the kink of the
$\ell_1$ norm at the origin. Hence, adapted optimization methods have to
be used, such as proximal splitting methods~\cite{Combettes2011}. In
PyHST, the iterative shrinkage-thresholding
algorithm~\cite{Daubechies2004} (ISTA) or the accelerated FISTA (fast ISTA) ~\cite{Beck2009}
algorithm can be used to solve the problem
(\ref{eq:opt_tv}).

In the current version of PyHST, the user can choose the value of $\beta$
and a number of ISTA or FISTA iterations. Convergence is typically
reached in a few hundreds of iterations, but FISTA converges faster than
ISTA because of the accelerated scheme. No automatic convergence
criterion is implemented yet, but the value of the energy in
(\ref{eq:opt_tv}) is displayed so that the user can check whether
satisfying convergence is reached or not. For the choice of $\beta$, a
higher value should be selected when the level of noise on the measures
is greater, but the optimal value of $\beta$ also depends on the number
of measurements, or of the magnitude of $\mathbf{x}$ (because the data
fidelity term and the total variation are not homogeneous in
$\mathbf{x}$), and on the sparsity of the image gradient, that is on its
microstructure. If one has access to the ground truth of representative
images from high-quality acquisitions, it is possible to calibrate the
optimal value $\beta$ in degraded acquisition conditions by maximizing
the signal over noise ratio of the reconstructed image. If the ground
truth is not available, statistical methods such as the discrepancy
principle~\cite{Wen2012} or generalized cross-validation~\cite{Liao2009}
can be used to select the optimal regularization parameter. Such methods
could be implemented in future versions of PyHST.

\subsection{Dictionary Learning}
For images non piece-wise constant the gradient is not sparse and therefore the total variation penalization is not the best choice.  On the opposite the intrinsic sparsity structure of a generic class of images can be learned  by the dictionary-learning technique~\cite{mairal}.
 This technique consists in building   an over-complete basis of functions, over  an $m \times m$  domain, such that, taken an $m \times m$ patch from an image  belonging to the studied image class, the patch can be approximated with good precision as a linear combination of a small number of  basis functions.  When this approximation is used we obtain a sparse representation for the images of that class.
  When a noisy image is represented  by the patch basis, the features of the original images will be accurately  fitted with a small number of components. The noise instead has in general no intrinsic sparsity, and if it happens to have a sparsity structure, it is with high probability very different from the sparsity structure of the original images. Therefore the noise will be reproduced only if we allow a large number of components (the patch basis is over-complete so it can represent the noise) but it will be effectively filtered out if we approximate the noisy image with a small number of components.

The patches denoising technique is often used with overlapping patches to avoid discontinuities at the patches borders. These discontinuities appear at the patches border because features which cross the patch eccentrically are weakly detected. A line crossing  the central region of a patch, as an example, will be detected if the basis of functions has been learned to fit such kind of features. 
A line crossing the patch in a corner point, instead, is equivalent to a noisy point.  The use of overlapping patches and the application of post-process averaging, with a weight  that depends on the distance from the patch center,  cures 
this problem and is widely used for image denoising producing good result~\cite{elad}.

In details, for each image patch, first, one optimizes an objective function which is composed of a fidelity term which forces the solution  image towards the data image and of a penalization term which promotes sparsity. The averaging is applied as post-processing once every patch has been fitted by  functional optimization procedures.

 In the case of tomography we want to address,  beside the denoising problem,  the  problem of reconstruction from a reduced amount of data. In this case, in fact, the use of a sparse representation of the reconstructed image can  fill the gap left by the missing data  with the  {\em a-priori} knowledge contained in the dictionary.

 A proper convex formulation including consistently the averaging over the overlapping zone,  has recently been introduced~\cite{overlap}.  
The convexity property is of paramount importance for the robustness of the numerical implementation.
 We  are going to remind here the basic principle of this method
and describe its GPU implementation in PyHST2.

\subsubsection{The Formalism}
  
We choose a set of patches which covers the whole image, and we allow for overlapping.
 In this case, for a given pixel, the sum of the patch indicator functions is greater or equal to one :

\begin{equation}
\sum_{p}{\bf 1}_{p}({\bf i})\geq1;\:\forall {\bf i}.
\end{equation}

where ${\bf 1}_{p}(i)$ is the indicator function for patch $p$. It is equal to $1$ if the pixel $\bf i$ (intended as a two dimensional vector)
belongs to the patch $p$ and is zero otherwise.

We denote by ${{\bf1}_{p}^{c}}$  the indicator functions for   patch core. The core is a  part of  the patch
which is selected  in order to  make a non-overlapping covering:
\begin{equation}
{{\bf 1}_{p}^{c}}(i)\leq{{\bf1}_{p}}(i);\quad\sum_p{{\bf 1}_{p}^{c}}(i)=1;\:\forall i.
\end{equation}

and, for a given point $i$, ${{\bf 1}_{p}^{c}}(i)$ indicates which patch
$p$ has its center $C_{p}$ closest to point $i$:

\begin{equation}
\sum_{p}{{\bf 1}_{p}^{c}}(i)\left\Vert i-C_{p}\right\Vert _{1}\leq\left\Vert i-C_{p^{\prime}}\right\Vert _{1};\quad\forall p^{\prime},i.
\end{equation}

The solution ${\bf x}$ is composed using the central part of the
patches as indicated by the functions ${{\bf 1}_{p}^{c}}$:

\begin{equation}
x_{i}=\sum_{p}{\bf 1}_{p}^{c}(i)\sum_{k}w_{kp}\varphi_{k}(i-r_{p}).
\end{equation}

Now we introduce the $\Pi$ operator which is the projection operator,
for tomography reconstruction, and is the identity for image denoising.
The functional $F({\bf w})$ whose minimum gives the optimal solution
is written, for both applications, as:

\begin{eqnarray}
F({\bf w})& =& f({\bf w})+g({\bf w});\quad g({\bf w})=\beta\left\Vert {\bf w}\right\Vert _{1};\nonumber\\
f({\bf w})& =& \left\Vert {\bf y}-{\bf \Pi}({\bf x})\right\Vert _{2}^{2}+\nonumber\\
& & \rho\sum_{pi}{\bf{1}_{p}(i)}\left(x_{i}-\sum_{k}w_{kp}\varphi_{k}(i-r_{p})\right)^{2}. ~~
\end{eqnarray}

 The factor $\rho$  weights a similarity-inducing term which
pushes all the overlapping patches, which touch a point $i$, toward
the value $x_{i}$ of the global solution $x$ in that point. 


The solution is found with the FISTA method, using the gradient of
$f({\bf w})$ which is easily written in compact form:

\begin{multline}
\label{gradient}
  \frac{\partial f({\bf w})}{\partial w_{kp}}=\sum_{i}2\varphi_{k}(i-r_{p})\bf{1}_{p}^{c}(i)\left\{ \left({\bf \Pi}^{T}\left({\bf \Pi}({\bf x})-{\bf y}\right)\right)_{i} 
 \begin{matrix} ~\\ ~ \\~ \end{matrix} \right. +\\
 \left.  \rho\sum_{p^{\prime}}\bf{1}_{p^{\prime}}(i)\left(x_{i}-\sum_{k^{\prime}}w_{k^{\prime}p^{\prime}}\varphi_{k^{\prime}}(i-r_{p^{\prime}})\right)\right\}+\\
  \sum_{i}2\varphi_{k}(i-r_{p})\rho{\bf 1}_{p}(i)\left(\sum_{k^{\prime}}w_{k^{\prime}p}\varphi_{k^{\prime}}(i-r_{p})-x_{i}\right)
\end{multline}

\subsubsection{The Implementation}

The ISTA and FISTA procedures consist in a gradient descend step followed by a shrinking step for the $L_1$ norm of the ${\bf w}$ coefficients .
While the $L_1$ norm shrinking step is trivial, the calculation of the gradient from equation \ref{gradient} 
is numerically expensive.

The most time expensive operations, for the gradient, are those involving the over-complete basis of functions. These operation have to be performed for every patch.
 There are two different operations involving the basis. One is the scalar product between a patch extracted from an image and all the basis functions for that patch. The other is the construction of an image patch from its coefficients on the basis.

 From the hardware point of view the most efficient way to perform these operations is by one single matrix-matrix multiplication where the concerned quantities are regrouped together. This operation can be efficiently performed with a single call to an optimized Basic Linear Algebra Subprograms (BLAS) library.
In our case we call the  Nvidia provided  cuBlas routines.
To do so we regroup patches in one single matrix ${\bf P}$ having dimensions $N \times m^2$, and containing all the $N$ extracted patches. The matrix for the basis functions (the dictionary) is denoted by the symbol  ${\bf D}$ and 
has dimension $N_c  \times m^2 $. It contains all the $N_c$ components of the over-complete basis (with $N_c \geq m^2$). The free coefficients are contained in a  single matrix 
${\bf W}$ having dimensions $N \times N_c$. The update of ${\bf W}$ is done a descent step along the gradient 
${\bf G}$ given by equation \ref{gradient} and implemented as shown in algorithm \ref{algogradient}

\label{aprioridl}

\section{Applications}
\label{application}
\begin{figure}
\caption{\label{emma} A reconstructed $2k \times 2k$ slice of a quartz grain sample. The reconstruction has been done
using a reduced set of only $150$ projection and a choice of three different $\beta$'s: from left to right
$\beta$ is equal to $10,3000$ and $1\times 10^5$. }
\begin{minipage}[t]{0.33\textwidth}%
\begin{center}
\includegraphics[width=0.95\columnwidth,height=0.95\columnwidth]{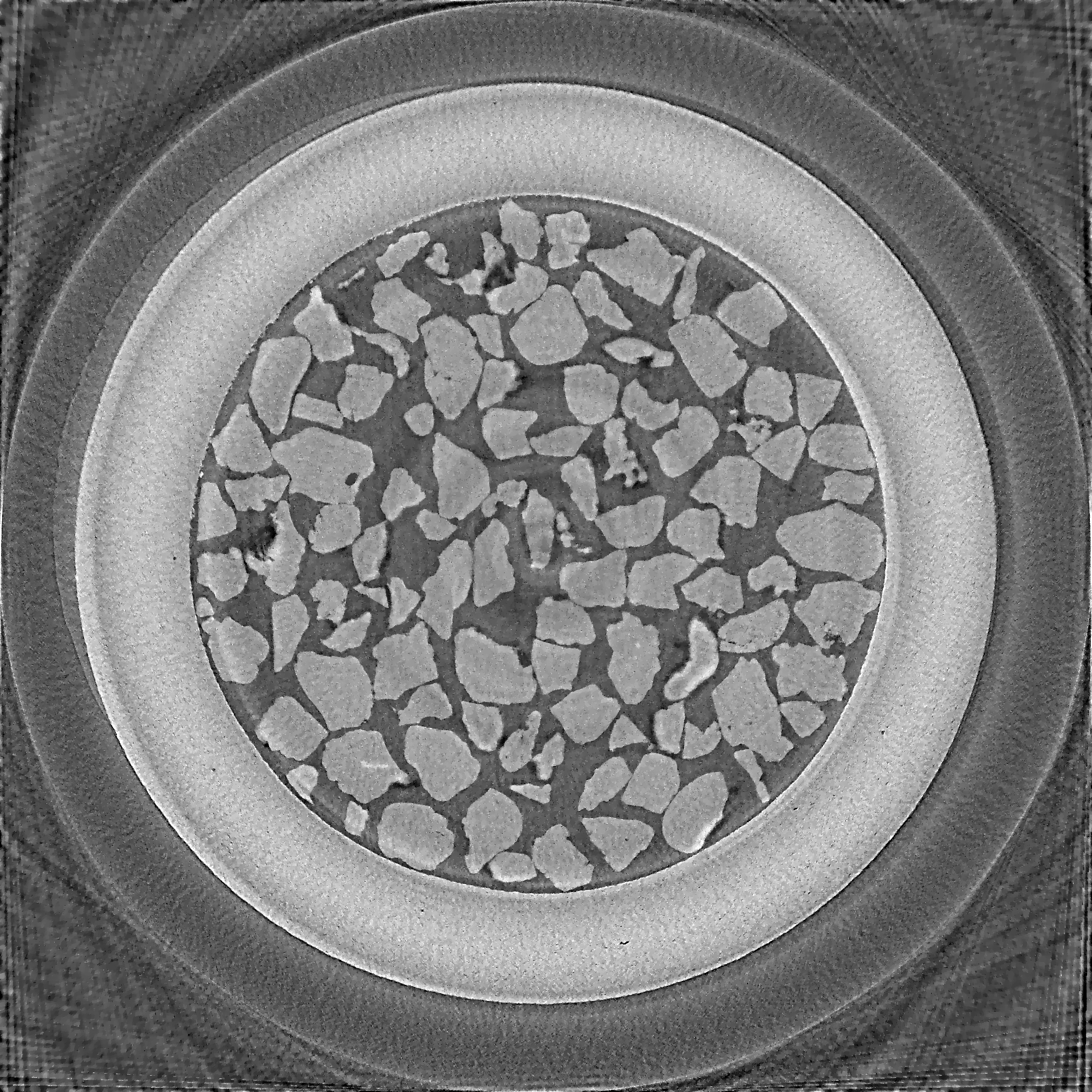}
\par\end{center}%
\end{minipage}
\hfill{}%
\begin{minipage}[t]{0.33\textwidth}%
\begin{center}
\includegraphics[width=0.95\columnwidth,height=0.95\columnwidth]{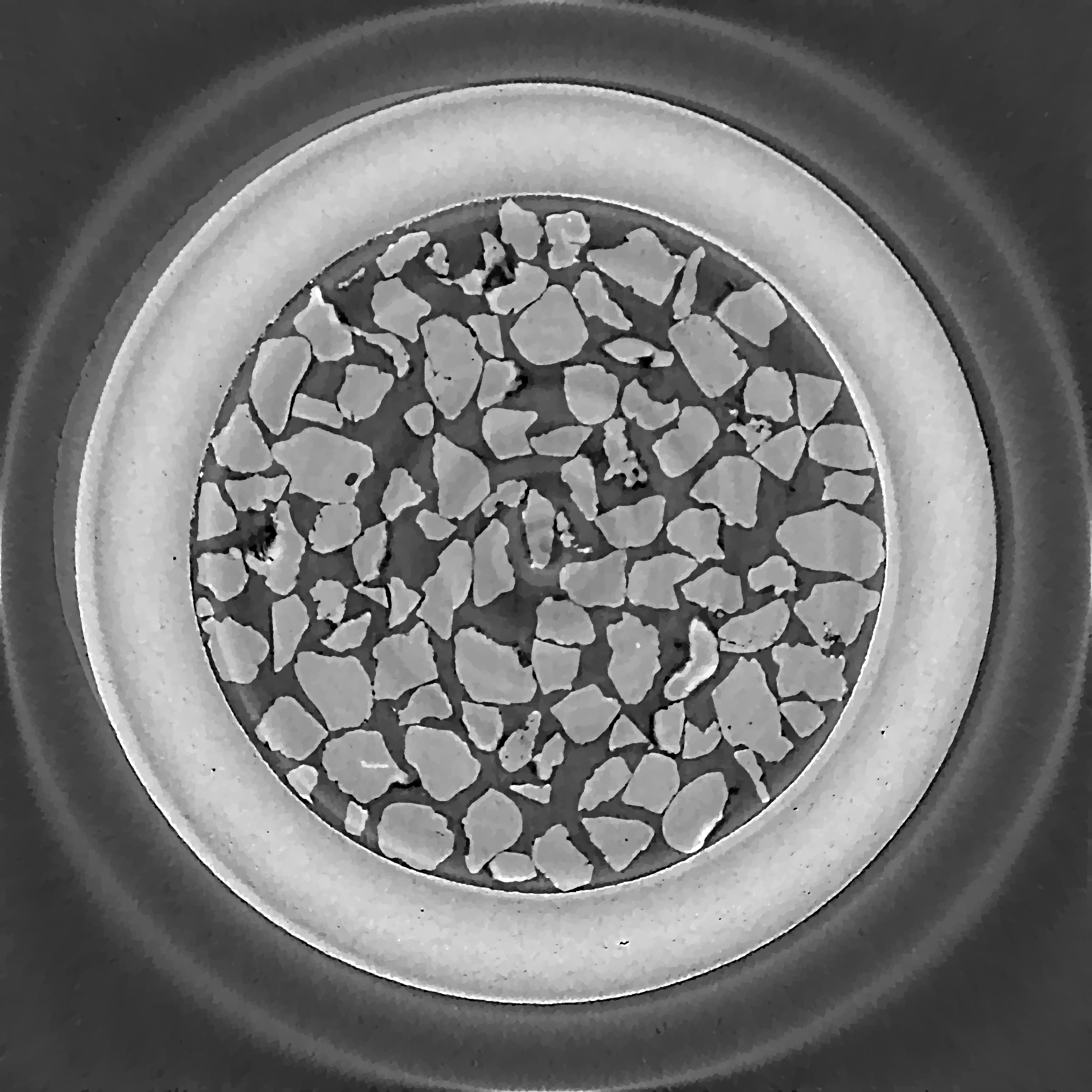}
\par\end{center}%
\end{minipage}\hfill{}%
\begin{minipage}[t]{0.33\textwidth}%
\begin{center}
\includegraphics[width=0.95\columnwidth,height=0.95\columnwidth]{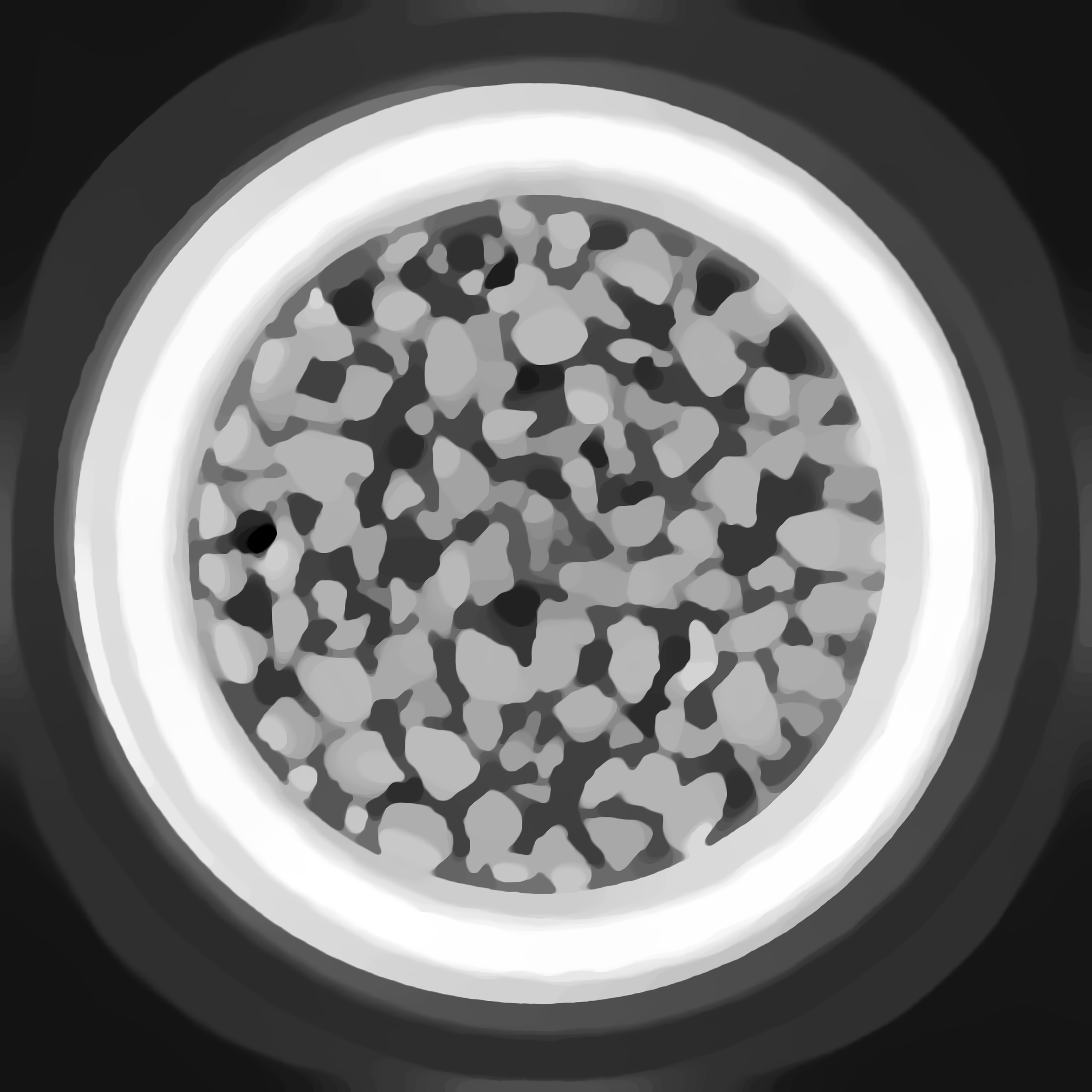}
\par\end{center}%
\end{minipage}
\end{figure}

\begin{figure}
\caption{\label{cross_cos} Plot, for the quartz grains sample, of two statistical estimators. The cross validation estimator is the distance from a selected
  projection which has not been used for the reconstruction. The decoherence maximising estimator is the cosinus between the removed noise and the obtained image.  }
\begin{minipage}[t]{0.5\textwidth}%
\begin{center}
\includegraphics[width=0.95\columnwidth,height=0.95\columnwidth]{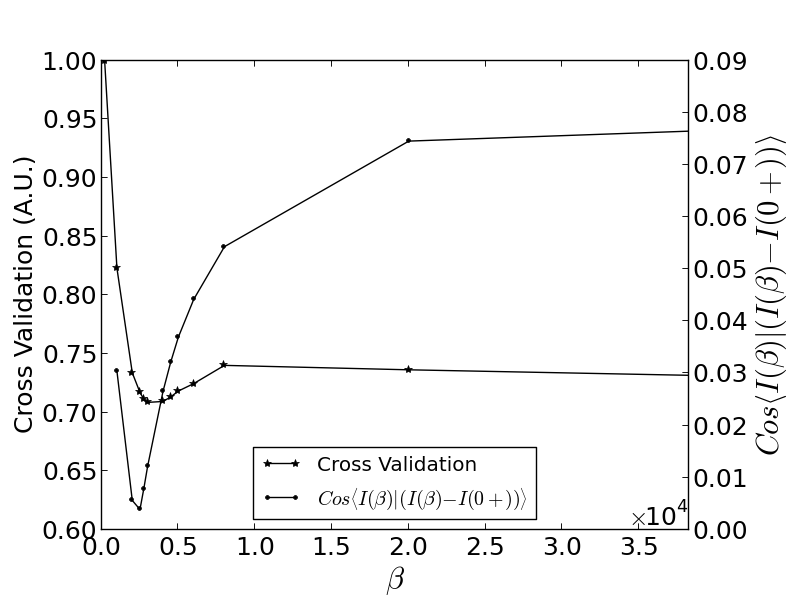}
\par\end{center}%
\end{minipage}
\hfill{}%
\end{figure}

We show in figure~\ref{emma} a reconstructed $2k \times 2k$ slice of a quartz grain sample. The reconstruction has been done
using a reduced set of only $150$ projection and a choice of three different $\beta$'s: from left to right
$\beta$ is equal to $10,3000$ and $1\times 10^5$. These values have been taken in order to encompass
the optimal $\beta$ that we can estimate from figure~\ref{cross_cos}.
In this figure we plot two quantities that we calculate in a postprocessing phase
after reconstructing the  sample at different $\beta$'s. One is the cross-validation estimator 
and the other is a new estimator that we define here and for which we
 coin the name {\it decoherence maximising estimator} $E_{cohe}$. 
The cross-validation estimator is obtained in the following canonical way. 
First we reconstruct the slice  using all the projections except a selected one at a given angle. Then
we calculate the quadratic distance  between  the reconstructed sample projection and the data for the excluded angle.

The  decoherence maximising estimator $E_{cohe}$  is constructed on the basis 
of the fact that the following cosinus is close to zero
\begin{equation}
    Cos\left( {\bf I}(\beta )-{\bf I}(0+), {\bf\hat{I}}\right)\simeq 0
\end{equation}
 for $\beta<\beta_{optimal}$ where ${\bf\hat{I} }$ is the unknown ground-truth, and where $0+$ is a very small value.
This affirmation relies on the fact that, if the penalization functional
has been properly chosen, at small  $\beta$'s  the penalization term removes a noisy signal
which is decoherent with the intrinsic sparsity structure . 

Because we don't know the ground truth we define our estimator $E_{cohe}$ as
\begin{equation}
  E_{cohe} = Cos\left( {\bf I}(\beta )-{\bf I}(0+), {\bf I}(\beta)\right)
\end{equation} 
A  thus defined estimator  has two asymptotic regions. On the left of $\beta_{optimal}$
$E_{cohe}$ must be far from zero because $( {\bf I}(\beta )-{\bf I}(0+))$ and ${\bf I}(\beta)$ both contain the noisy signal.
On the right, instead, $( {\bf I}(\beta )-{\bf I}(0+))$ contains part of the optimal solution because 
at strong $\beta$ the penalisation induces distortion.
We argue therefore that the  minimum of $abs(E_{cohe})$ must lie not far from the optimal $\beta$.
We see from figure~\ref{cross_cos} that the two estimators minima are not far from each other.
The value $\beta=3000$ used for the central image of figure~\ref{emma} corresponds to the minimum
of the cross-validation estimator. The two side images illustrate the image distortion
at high $\beta$ on the right at $\beta=1\times 10^5$, and the noise removal going
from $\beta=10$ to the optimal image.  

\begin{figure}
\caption{\label{flower} Reconstruction of a $2k\times 2k$  using 150 projections  with added gaussian white noise.
 The $\beta$ values  are from left to right $\beta=0.001$,$\beta=0.065$  corresponding to the ground-truth minimal distance,
 and $\beta=0.035$  corresponding to  the minimum of the decoherence maximising estimator.}
\begin{minipage}[t]{0.3\textwidth}%
$a) \beta=0.001$
\begin{center}
\includegraphics[width=0.95\columnwidth,height=0.95\columnwidth]{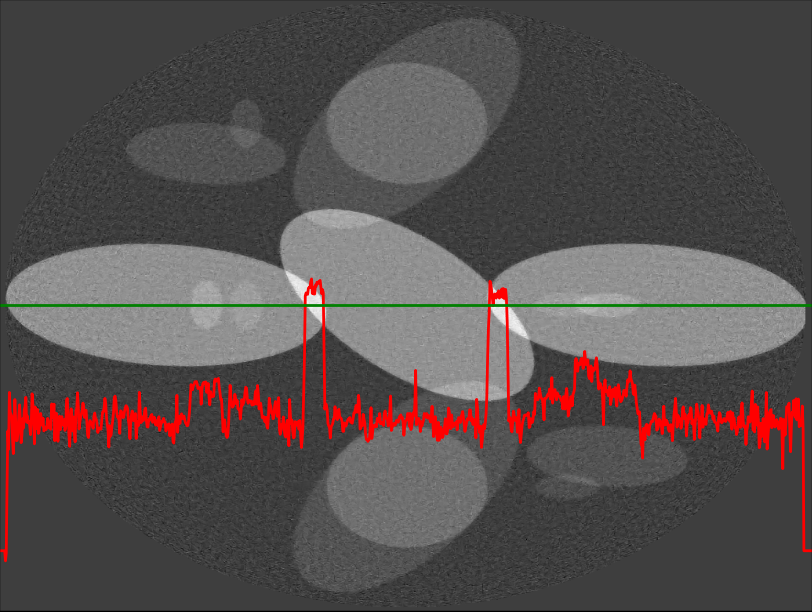}
\par\end{center}%
\end{minipage}\hfill{}%
\hfill{}%
\begin{minipage}[t]{0.3\textwidth}%
$b) \beta=0.065$
\begin{center}
\includegraphics[width=0.95\columnwidth,height=0.95\columnwidth]{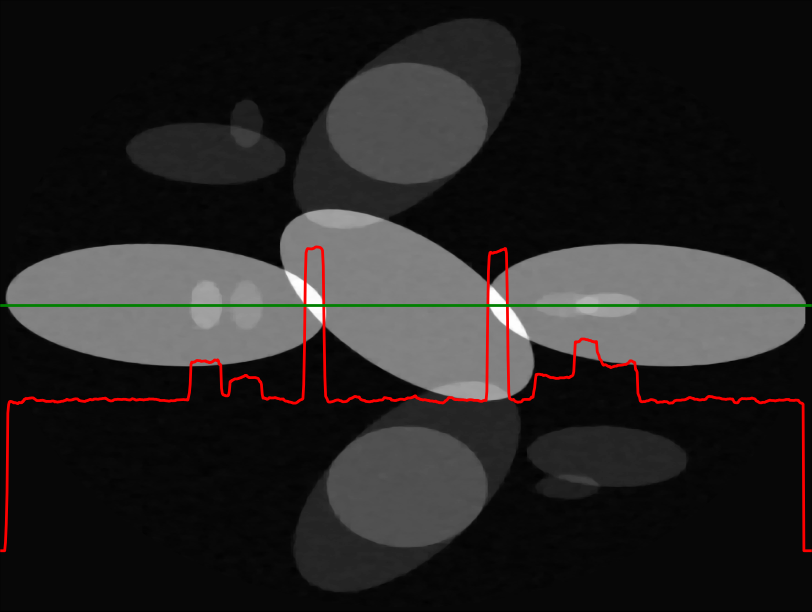}
\par\end{center}%
\end{minipage}
\hfill{}%
\begin{minipage}[t]{0.3\textwidth}%
$c) \beta=0.035$
\begin{center}
\includegraphics[width=0.95\columnwidth,height=0.95\columnwidth]{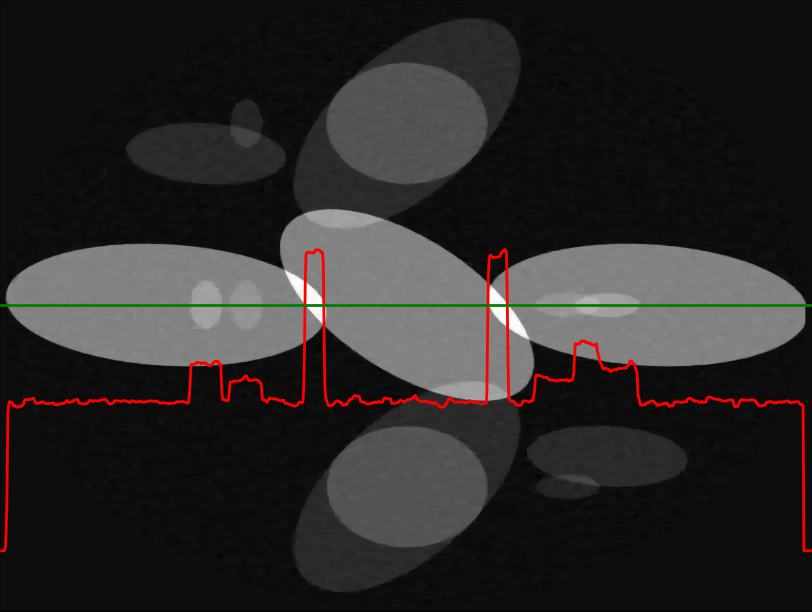}
\par\end{center}%
\end{minipage}
\end{figure}

\begin{figure}
\begin{minipage}[t]{0.4\textwidth}%
\begin{center}
\caption{\label{patches}The basis of patches for the phantom reconstruction}
\includegraphics[width=0.95\columnwidth,height=0.95\columnwidth]{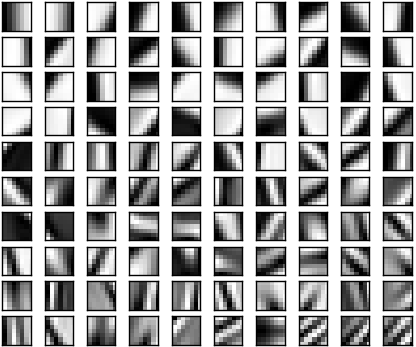}
\par\end{center}%
\vskip-1cm
\end{minipage}
\end{figure}

\begin{figure}
\begin{minipage}[t]{0.5\textwidth}%
\begin{center}
\caption{\label{cross_cos_truth} Plot of the statistical estimators and of the distance from ground-truth.
 The estimators give a  $\beta$ which is smaller than, but still close to the optimal one. }
\includegraphics[width=0.95\columnwidth,height=0.95\columnwidth]{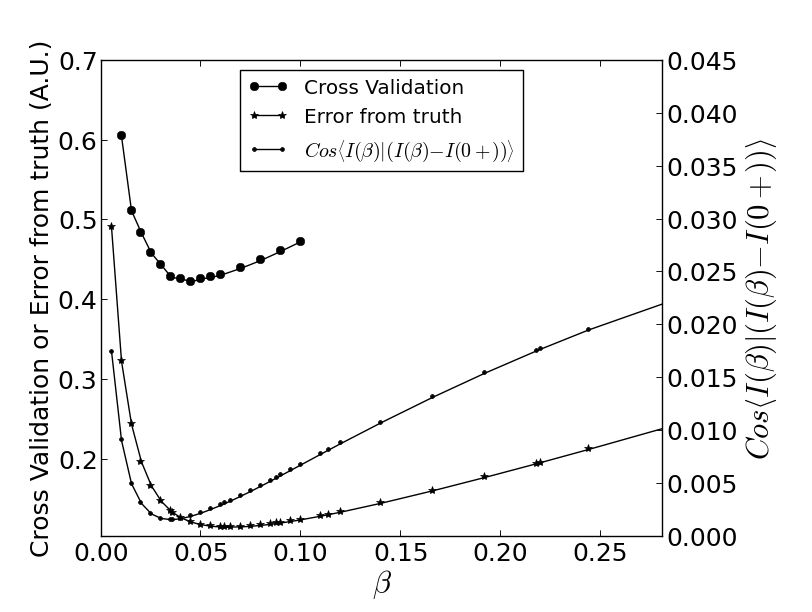}
\par\end{center}%
\end{minipage}
\end{figure}

For this experimental sample we don't know the ground truth. In order to validate
the use of the two estimators we reconstruct a $2k\times 2k$   phantom, whose reconstruction at different $\beta$'s
is shown in figure~\ref{flower}. We apply the overlapping patches functional 
to provide at the same time an illustration of this new method.
We use 150 projections of a synthesised sinogram with added Gaussian white noise.
 The $\beta$ used values  correspond
from left to right to $0+$ (calculation done at $\beta=0.001$), to the ground-truth minimal distance
 at $\beta=0.065$ and to the minimum of the maximal decoherence estimator
 at $\beta=0.035$.
 The basis of patches is shown in figure~\ref{patches}. The plot of the estimators and of the ground-truth distance
is shown in figure~\ref{cross_cos_truth}, where we have varied $\beta$ while keeping $\rho$ fixed and using 
the same basis of patches as in~\cite{overlap}, shown in figure~\ref{patches}.
 We can see that the estimators minima are both not far
from the ground-truth optimal value, and are close to each other. The error done
using the estimator can be checked on image~\ref{flower}b and ~\ref{flower}c.
The decrease in image quality, between the optimal to suboptimal values, is barely detectable by the eye.

\section{Conclusions}
\label{conclusions}
 We have presented the PyHST code which is now distributed under GPL license~\cite{pyhst2_gpl,GPL}.
This code implements advanced {\em a-priori} knowledge techniques and we have discussed their implementations.
For the new, recently found, dictionary-learning technique based on overlapping patches we have
provided the scheme for its efficient GPU implementation.
We have tested the {\em a-priori} techniques and we have applied statistical methods to validate the 
parameters choices. In particular we have applied, beside cross-validation, 
an original {\em a-priori} validation technique based on decoherence maximisation. 



\bibliographystyle{elsarticle-num}
\bibliography{bib}

\begin{thebibliography}{10}
\expandafter\ifx\csname url\endcsname\relax
  \def\url#1{\texttt{#1}}\fi
\expandafter\ifx\csname urlprefix\endcsname\relax\def\urlprefix{URL }\fi
\expandafter\ifx\csname href\endcsname\relax
  \def\href#1#2{#2} \def\path#1{#1}\fi

\bibitem{Brunet2004277}
M.~Brunet, F.~Guy, J.-R. Boisserie, A.~Djimdoumalbaye, T.~Lehmann, F.~Lihoreau,
  A.~Louchart, M.~Schuster, P.~Tafforeau, A.~Likius, H.~T. Mackaye, C.~Blondel,
  H.~Bocherens, L.~D. Bonis, Y.~Coppens, C.~Denis, P.~Duringer, V.~Eisenmann,
  A.~Flisch, D.~Geraads, N.~Lopez-Martinez, O.~Otero, P.~P. Campomanes,
  D.~Pilbeam, M.~P. de~León, P.~Vignaud, L.~Viriot, C.~Zollikofer, Toumai,
  mioc\`ene sup\'erieur du tchad, le nouveau doyen du rameau humain, Comptes
  Rendus Palevol 3~(4) (2004) 277 -- 285.
\newblock \href {http://dx.doi.org/10.1016/j.crpv.2004.04.004}
  {\path{doi:10.1016/j.crpv.2004.04.004}}.

\bibitem{Salvo2003}
L.~Salvo, P.~Cloetens, E.~Maire, S.~Zabler, J.~Blandin, J.-Y. Buffi{\`e}re,
  W.~Ludwig, E.~Boller, D.~Bellet, C.~Josserond, X-ray micro-tomography an
  attractive characterisation technique in materials science, Nuclear
  instruments and methods in physics research section B: Beam interactions with
  materials and atoms 200 (2003) 273--286.

\bibitem{Zhao2012}
Y.~Zhao, E.~Brun, P.~Coan, Z.~Huang, A.~Sztr\'{o}kay, P.~C. Diemoz,
  S.~Liebhardt, A.~Mittone, S.~Gasilov, J.~Miao, A.~Bravin,
  \href{http://www.pnas.org/cgi/content/long/1204460109v1}{{High-resolution,
  low-dose phase contrast X-ray tomography for 3D diagnosis of human breast
  cancers.}}, Proceedings of the National Academy of Sciences of the United
  States of America 109~(45) (2012) 18290--18294.
\newblock \href {http://dx.doi.org/10.1073/pnas.1204460109}
  {\path{doi:10.1073/pnas.1204460109}}.
\newline\urlprefix\url{http://www.pnas.org/cgi/content/long/1204460109v1}

\bibitem{overlap}
A.~Mirone, E.~Brun, P.~Coan, {A Convex Functional for image denoising based on
  Patches with Constrained Overlaps and its vectorial application to Low Dose
  Differential Phase Tomography}, submitted to IEEE proceeding on image
  processing vvv (yyy) xxx--xxx.
\newblock \href {http://arxiv.org/abs/1305.1256} {\path{arXiv:1305.1256}}.

\bibitem{nano}
P.~Bleuet, P.~Cloetens, P.~Gergaud, D.~Mariolle, N.~Chevalier, R.~Tucoulou,
  J.~Susini, A.~Chabli, \href{http://dx.doi.org/10.1063/1.3117489}{{A hard
  x-ray nanoprobe for scanning and projection nanotomography}}, Review of
  Scientific Instruments 80~(5) (2009) 056101+.
\newblock \href {http://dx.doi.org/10.1063/1.3117489}
  {\path{doi:10.1063/1.3117489}}.
\newline\urlprefix\url{http://dx.doi.org/10.1063/1.3117489}

\bibitem{paganin}
D.~Paganin, S.~C. Mayo, T.~E. Gureyev, P.~R. Miller, S.~W. Wilkins,
  \href{http://dx.doi.org/10.1046/j.1365-2818.2002.01010.x}{Simultaneous phase
  and amplitude extraction from a single defocused image of a homogeneous
  object}, Journal of Microscopy 206~(1) (2002) 33--40.
\newblock \href {http://dx.doi.org/10.1046/j.1365-2818.2002.01010.x}
  {\path{doi:10.1046/j.1365-2818.2002.01010.x}}.
\newline\urlprefix\url{http://dx.doi.org/10.1046/j.1365-2818.2002.01010.x}

\bibitem{Candes2006}
E.~Cand{\`e}s, J.~Romberg, T.~Tao, Robust uncertainty principles: Exact signal
  reconstruction from highly incomplete frequency information, Information
  Theory, IEEE Transactions on 52~(2) (2006) 489--509.

\bibitem{Candes2006b}
E.~Candes, J.~Romberg, T.~Tao, Stable signal recovery from incomplete and
  inaccurate measurements, Communications on pure and applied mathematics
  59~(8) (2006) 1207--1223.

\bibitem{Song2007}
J.~Song, Q.~Liu, G.~Johnson, C.~Badea, {Sparseness prior based iterative image
  reconstruction for retrospectively gated cardiac micro-CT}, Medical physics
  34~(11) (2007) 4476.

\bibitem{Sidky2008}
E.~Sidky, X.~Pan, Image reconstruction in circular cone-beam computed
  tomography by constrained, total-variation minimization, Physics in medicine
  and biology 53 (2008) 4777.

\bibitem{Herman2008}
G.~Herman, R.~Davidi, Image reconstruction from a small number of projections,
  Inverse Problems 24 (2008) 045011.

\bibitem{Tang2009}
J.~Tang, B.~Nett, G.~Chen, Performance comparison between total variation
  (tv)-based compressed sensing and statistical iterative reconstruction
  algorithms, Physics in Medicine and Biology 54 (2009) 5781.

\bibitem{Sidky2010}
E.~Sidky, M.~Anastasio, X.~Pan, {Image reconstruction exploiting object
  sparsity in boundary-enhanced X-ray phase-contrast tomography}, Optics
  Express 18~(10) (2010) 10404--10422.

\bibitem{Jia2010}
X.~Jia, Y.~Lou, R.~Li, W.~Song, S.~Jiang, Gpu-based fast cone beam ct
  reconstruction from undersampled and noisy projection data via total
  variation, Medical physics 37 (2010) 1757.

\bibitem{Chambolle2004}
A.~Chambolle, An algorithm for total variation minimization and applications,
  Journal of Mathematical imaging and vision 20~(1-2) (2004) 89--97.

\bibitem{Combettes2011}
P.~Combettes, J.~Pesquet, Proximal splitting methods in signal processing,
  Fixed-Point Algorithms for Inverse Problems in Science and Engineering (2011)
  185--212.

\bibitem{Daubechies2004}
I.~Daubechies, M.~Defrise, C.~De~Mol, An iterative thresholding algorithm for
  linear inverse problems with a sparsity constraint, Communications on pure
  and applied mathematics 57~(11) (2004) 1413--1457.

\bibitem{Beck2009}
A.~Beck, M.~Teboulle, Fast gradient-based algorithms for constrained total
  variation image denoising and deblurring problems, Image Processing, IEEE
  Transactions on 18~(11) (2009) 2419--2434.

\bibitem{Wen2012}
Y.-W. Wen, R.~H. Chan, Parameter selection for total-variation-based image
  restoration using discrepancy principle, Image Processing, IEEE Transactions
  on 21~(4) (2012) 1770--1781.

\bibitem{Liao2009}
H.~Liao, F.~Li, M.~K. Ng, Selection of regularization parameter in total
  variation image restoration, JOSA A 26~(11) (2009) 2311--2320.

\bibitem{mairal}
J.~Mairal, F.~Bach, J.~Ponce, G.~Sapiro, Online learning for matrix
  factorization and sparse coding, Journal of Machine Learning Research 11
  (2010) 19--60.

\bibitem{elad}
M.~Elad, M.~Aharon, \href{http://dx.doi.org/10.1109/tip.2006.881969}{Image
  denoising via sparse and redundant representations over learned
  dictionaries}, Image Processing, IEEE Transactions on 15~(12) (2006)
  3736--3745.
\newblock \href {http://dx.doi.org/10.1109/tip.2006.881969}
  {\path{doi:10.1109/tip.2006.881969}}.
\newline\urlprefix\url{http://dx.doi.org/10.1109/tip.2006.881969}

\bibitem{pyhst2_gpl}
A.~Mirone, \href{http://forge.epn-campus.eu/projects/pyhst2}{Gpl release of
  pyhst2.}, \url{http://forge.epn-campus.eu/projects/pyhst2} (2013).
\newline\urlprefix\url{http://forge.epn-campus.eu/projects/pyhst2}

\bibitem{GPL}
Gnu general public license, version 2,
  \url{http://www.gnu.org/licenses/gpl-2.0.html}.

\end{thebibliography}







\begin{algorithm}
\caption {Gradient computation.  $\bf D$ 
has dimension $N_c  \times m^2 $ and contains all the $N_c$ components of the over-complete basis ( with $N_c \geq m^2$ ).  $\bf P, P^\prime$ have dimensions $N \times m^2$, and contain a whole set of  $N$ overlapping patches  each.  $\bf W$ has dimensions $N \times N_c$ and contains the coefficients and $\bf G$ is the gradient of equation \ref{gradient}.}\label{algogradient}

\begin{algorithmic}
\State ${\bf P } \leftarrow  {\bf W} \cdot {\bf D}$;

\ForAll { image positions {\bf j}}
\State $ p \leftarrow C_{\bf j}$; \Comment{{\bf C} is precalculated and  $C_{\bf j}$  is the patch whose core contains $\bf j$ }
\State ${\bf i} \leftarrow {\bf j}-{\bf r}_{ p}$;
\State  $X_{\bf j} \leftarrow P_{p{\bf i}}$ ; \Comment{Build the solution image $\bf X$ from patches cores}
\EndFor
\State ${\bf R}= {\widehat {\bf Bp}}(  data- {\widehat {\bf Fp}}({\bf X})  )$ ; \Comment $\bf R$ is the backprojection of the residual error
\State \Comment{${\widehat {\bf Bp}}$ and ${\widehat {\bf Fp}}$ are the back-projection and projection operators.}
\State ${\bf S} \leftarrow 0$; \Comment{$\bf S$ is a temporary image for the patches sum}

\ForAll {patches $p$ and intra-patch positions ${\bf i}$}
\State $ S_{{\bf r}_p+{\bf i}} \leftarrow  S_{{\bf r}_p+{\bf i}} + P_{p \bf i}$;  \Comment{Build the sum of all the patches covering ${\bf i}$}
\EndFor
\State ${\bf P^\prime} \leftarrow 0$; \Comment{${\bf P}^\prime$ is an auxiliary array of patches}
\ForAll { $p$ ,${\bf i}$}    \Comment{Build the patch array P from the residue R}
\If{ $p == C_{r_p+i}$}  
\State \Comment{$N_v$ is the number of patches touching ${\bf i}$}
\State $P^\prime_{p \bf i} \leftarrow P^\prime_{p \bf i}+R_{{\bf r}_p+{\bf i}}+ \rho(   N_v * X_{{\bf r}_p+{\bf i}} -S_{{\bf r}_p+{\bf i}} )$ ;  
\EndIf  
\State  $P^\prime_{p \bf i} \leftarrow  P^\prime_{p \bf i} +  \rho(  P_{{\bf r}_p+{\bf i}} -  X_{{\bf r}_p+{\bf i}} )$;
\EndFor
\State ${\bf G} \leftarrow  2 {\bf P^\prime} \cdot {\bf D}^T$;
\end{algorithmic}
\end{algorithm}

\newpage


\end{document}